\input amstex
\magnification=1200
\documentstyle{amsppt}
\NoRunningHeads
\NoBlackBoxes
\define\End{\operatorname{End}}
\define\<{\left<}
\define\>{\right>}
\define\tr{\operatorname{tr}}
\define\FR{\operatorname{\Cal F\Cal R}}
\define\df{\operatorname{def}}
\define\sltwo{\operatorname{\frak s\frak l}(2,\Bbb C)}
\define\asltwo{\operatorname{\frak a\frak s\frak l}(2,\Bbb C)}
\define\Mat{\operatorname{Mat}}
\topmatter
\title
Cross-projective representations of pairs of anticommutative
algebras, alloys and finite-dimensional irreducible representations of
some infinite-dimensional Lie algebras
\endtitle
\author Denis V.Juriev
\endauthor
\date math.RT/9806005\qquad\qquad\qquad June 02, 1998\enddate
\endtopmatter
\document
This article is devoted to some ``strange'' phenomena in representation
theory, whose existence seems to be unexplainable to me. The discussed
subjects are, perphaps, rather far from the wide roads of classical and
modern researches, how\-ever, one might pay an attention to them as for
themselves as for any associations, which they may provoke. Concerning
the applications of the collected facts I as a realist suspect that
everything being correctly imagined really exists somewhere and if any
Theory is beyond our Practice one should follow the most practical advice
that is proposed, namely, let try to make the least a bit wider and do not
cut off the first.

\subhead 1.1. Cross-projective representations of pairs of anticommutative
al\-geb\-ras\endsubhead

\definition{Definition 1} Let $\frak A_i$ ($i=1,2$) be two anticommutative
algebras with brackets $\<\cdot,\cdot\>_i$. {\it A cross-projective
representation\/} of the pair $(\frak A_1,\frak A_2)$ in the linear space
$V$ is the pair $(T_1,T_2)$ of mappings $T_i:\frak A_i\mapsto\End(V)$ such
that
$$\aligned
[T_1(X),T_1(Y)]-T_1(\<X,Y\>_1)\equiv0&\mod T_2(\frak A_2),\\
[T_2(A),T_2(B)]-T_2(\<A,B\>_2)\equiv0&\mod T_1(\frak A_1),
\endaligned
$$
for all $X,Y\in\frak A_1$, $A,B\in\frak A_2$. Here $[\cdot,\cdot]$ is
the commutator of operators.
\enddefinition

The definition 1 is related to one of $\frak A$-projective representations
(see e.g.[1,2]). Any $\frak A$-projective representation of an anticommutative
algebra $\frak g$ defines a cross-projective representation of the pair
$(\frak g,\frak A_{[\cdot,\cdot]})$ (here $\frak A_{[\cdot,\cdot]}$ is the
commutator Lie algebra of the associative algebra $\frak A$).

\remark{Example 1} Let $\frak g$ be a Lie algebra decomposed into the
direct sum $\frak g_1\oplus\frak g_2$ of linear spaces supplied by the
brackets $\<\cdot,\cdot\>_i$, which are projections of the commutator
$[\cdot,\cdot]$ in $\frak g$ onto $\frak g_i$ along the complementary
subspace. Then any representation $T$ of $\frak g$ defines the
cross-projective representation $(T_1,T_2)$ ($T_i=\left.T\right|_{\frak
g_i}$) of the pair $(\frak g_1,\frak g_2)$.
\endremark

\remark{Remark 1} Definition 1 admits natural higher combinatorial
generalizations if one use a set of arbitrary number of anticommutative
algebras instead of their pair.
\endremark

One may introduce the basic concepts of representation theory (reducible,
ir\-re\-du\-cible, decomposable, etc. representations) for the cross-projective
representations. Sometimes the problem of classification of irreducible
finite-dimensional rep\-re\-sen\-tations of a pair of anticommutative algebras
is hugely wild, e.g. for a pair of abelian algebras $(\Bbb C^n,\Bbb C^m)$.
However, we shall see below that this problem is slightly timer for more
nontrivial pairs.

\remark{Remark 2} Let $(T_1,T_2)$ be a finite-dimensional representation of
the pair $(\frak g_1,\frak g_2)$ in the space $V$ then
$(T^\times_1,T^\times_2)$ ($T^\times_i(X)=T_i(X)-\dfrac{\tr T_i(X)}{\dim V}E$,
$E$ is the unit matrix) is also a representation of this pair.
\endremark

\remark{Remark 3} Let $\frak g$ be an anticommutative algebra then
$\bigwedge^2(\frak g)$ also possesses a natural structure of anticommutative
algebra so that it is possible to associate cross-pro\-jec\-tive
representations of the pair $(\frak g,\bigwedge^2(\frak g))$ with any
anticommutative algebra $\frak g$.
\endremark

Note that cross-projective representations do not form neither tensor nor
even abelian category. To improve the situation one should consider it from
a different point of view.

\subhead 1.2. Alloys and their representations\endsubhead

\definition{Definition 2}

{\bf A.} Let $\frak g$ be a linear space decomposed into the sum $\frak g_1
\biguplus\ldots\biguplus \frak g_n$ of its subspaces $\frak g_i$ (the sum is
not direct in general). Let us consider the subspace $W$ of $\bigwedge^2(\frak
g)$ of the form $\bigwedge^2(\frak g_1)\biguplus\ldots\biguplus\bigwedge^2
(\frak g_n)$. The linear space $\frak g$ supplied by the partial
anticommutative binary operation $[\cdot,\cdot]:W\mapsto \frak g$ is called
{\it an alloy}.

{\bf B.} {\it A representation\/} $T$ of the alloy $(\frak g,[\cdot,\cdot])$
in the linear space $V$ is the mapping $T:\frak g\mapsto\End(V)$ such that
$[T(X),T(Y)]=T([X,Y])$ for all pairs $(X,Y)$ from $W$. Here $[\cdot,\cdot]$
is the commutator of operators in the right hand side of the identity.
\enddefinition

\remark{Remark 4} The definition may be naturally generalized to arbitrary
subspaces $W$ of $\bigwedge^2(V)$ instead of ones of a specific form considered
above.
\endremark

Note that each Lie composite (see e.g.[1,3]) may be considered as an alloy.

\proclaim{Theorem 1} Let $(T_1,T_2)$ be a strict cross-projective
representation of the pair of anticommutative algebras $(\frak g_1,\frak g_2)$
such that $T_1(\frak g_1)\cap T_2(\frak g_2)=0$. Then $T=T_1\oplus T_2$
realizes a strict representation of a uniquely defined alloy
$\frak g=\frak g_1\oplus\frak g_2$ such that projections of the binary
operation $[\cdot,\cdot]$ in $\frak g$ on $\frak g_i$ along the complementary
subspaces coincide with the binary operations in the anticommutative algebras
$\frak g_i$.
\endproclaim

\remark{Remark 5} If one omits the condition of transversality of images of
$T_1$ and $T_2$ or the condition of strictness of representations then it
will be possible to construct an alloy but the least will not be uniquely
defined.
\endremark

\proclaim{Proposition 1} The representations of a fixed alloy $\frak g$ form
a tensor abelian ca\-te\-go\-ry.
\endproclaim

\subhead 1.3. Quaternary algebras and their alloyability
\endsubhead

\proclaim{Proposition 2} Let $\frak g=\frak g_1\oplus\frak g_2$ be an alloy
then $\frak g_i$ are supplied by the structure of both binary and
quaternary algebras. The quaternary operation $\<\cdot,\cdot,\cdot,\cdot\>_i:
\bigwedge^2(\bigwedge^2(\frak g_i))\mapsto\frak g_i$ is constructed in the
following way
$$\aligned
\<U,V,X,Y\>_1&=\lambda_2(\lambda_1(U,V),\lambda_1(X,Y)),\\
\<A,B,C,D\>_2&=\lambda_1(\lambda_2(A,B),\lambda_2(C,D)),
\endaligned$$
where $U,V,X,Y\in\frak g_1$, $A,B,C,D\in\frak g_2$ and $\lambda_i$ are
projections of the operation $[\cdot,\cdot]$ in the alloy $\frak g$
restricted to $\frak g_i$ and projected along $\frak g_i$
onto the complementary subspace.
\endproclaim

\proclaim{Corollary} Any strict cross-projective representation $(T_1,T_2)$
of the pair of an\-ti\-com\-mu\-ta\-tive algebras $(\frak g_1,\frak g_2)$
such that $T_1(\frak g_1)\cap T_2(\frak g_2)=0$ supplies both algebras by
additional quaternary operations.
\endproclaim

\remark{Remark 6} In general, the quaternary operations are not correlated
with binary operations in any way.
\endremark

Corollary is a manifestation of an interesting phenomena that a relation
of two objects (anticommutative algebras here) in their mutual
representation (cross-projective representation here) supplies them by
additional abstract algebraic struc\-tures (quaternary operations here).

\remark{Remark 7} Note that for any quaternary algebra $\frak g$ the following
statements hold:
\roster
\item"--" $\bigwedge^2(\frak g)$ is a quaternary algebra;
\item"--" $\bigwedge^2(\bigwedge^2(\frak g))$ is an extension of $\frak g$,
i.e. there exists a homomorphism $\pi:\bigwedge^2(\bigwedge^2(\frak g))\mapsto
\frak g$ of quaternary algebras (defined just by the quaternary operation
$\<\cdot,\cdot,\cdot,\cdot\>$ in $\frak g$).
\endroster
\endremark

\definition{Definition 3} Two quaternary algebras $\frak g_i$ ($i=1,2$) with
operations $\<\cdot,\cdot,\cdot,\cdot\>_i:\bigwedge^2(\bigwedge^2(\frak
g_i))\mapsto\frak g_i$ will be called {\it alloyable\/} iff they possess a
mutual factorization
$$\aligned
\<U,V,X,Y\>_1&=\lambda_2(\lambda_1(U,V),\lambda_1(X,Y)),\\
\<A,B,C,D\>_2&=\lambda_1(\lambda_2(A,B),\lambda_2(C,D)),
\endaligned$$
where $U,V,X,Y\in\frak g_1$, $A,B,C,D\in\frak g_2$, $\lambda_1:\bigwedge^2
(\frak g_1)\mapsto\frak g_2$ and $\lambda_2:\bigwedge^2(\frak g_2)\mapsto
\frak g_1$.
\enddefinition

\remark{Remark 8} The problem of an alloyability of two quaternary
algebras may be considered as a certain specific algebraic counterpart of
the Kolmogorov-Arnold problem for functions of several variables [4,5].
\endremark

\proclaim{Proposition 3} For any quaternary algebra $\frak g_1$ there
exists an alloyable quaternary algebra $\frak g_2$.
\endproclaim

\demo{Proof} One may consider $\bigwedge^2(\frak g_1)$ as $\frak g_2$.
\enddemo

The proposition 3 means that the relation of alloyability $R_A$ on the class
$\Cal Q$ of quaternary algebras in nondegenerate. Note that the dimension of
the space $\Cal Q_n$ of quaternary algebras of dimension $n$ is equal to
$\frac18(n+1)n^2(n-1)(n-2)$ whereas the dimension of the graph of the
relation $R_A$ in $\Cal Q_n\oplus\Cal Q_m$ is equal to $\frac12nm(n+m-2)$
so the condition of alloyability is rather strong. I suspect that the
problem of classification of pairs of alloyable quaternary algebras is wild.

\remark{Remark 9} Let $\frak g$ be an anticommutative algebra with the
bracket $\<\cdot,\cdot\>$ supplied by the additional structure of a quaternary
algebra with the operation $\<\cdot,\cdot,\cdot,\cdot\>$ then
$\bigwedge^2(\frak g)$ is also an anticommutative algebra supplied by such
additional structure (remarks 3,7). The quaternary operations in
$\frak g$ and $\bigwedge^2(\frak g)$ are alloyable, their fac\-to\-ri\-zations
$\lambda_1$, $\lambda_2$ coupled with the anticommutative binary operations
in $\frak g$ and $\bigwedge^2(\frak g)$ define the structure of an alloy in
the space $\frak g\oplus\bigwedge^2(\frak g)$.
\endremark

Note that in any anticommutative algebra $\frak g$ with the binary
operation $\<\cdot,\cdot\>$ one may construct an additional quaternary
operation $\<\cdot,\cdot,\cdot,\cdot\>$ as
$$\<A,B,C,D\>=\<\<A,B\>,\<C,D\>\>.$$
If the space $\frak g$ possesses two brackets $\<\cdot,\cdot\>_i$
($i=1,2$) then one may construct six quaternary operations
$\<\cdot,\cdot,\cdot,\cdot\>_{ij}^k$ ($i,j,k=1,2$) as
$$\<A,B,C,D\>_{ij}^k=\<\right.\<A,B\>_i,\<C,D\>_j\left.\>_k+
(i\leftrightarrow j).$$

\subhead 1.4. Universal envelopping Lie algebras of alloys and their
rep\-re\-sen\-tations
\endsubhead
Let $\frak g=\frak g_1\biguplus\ldots\biguplus\frak g_n$ be an alloy.
The Lie algebra $\tilde\frak g$ will be called the envelopping Lie algebra
for the alloy $\frak g$ iff there exists a monomorphism of $\frak g$
into $\tilde\frak g$ (the least is considered as an alloy) in the category
of alloys. The universal object in the category of envelopping Lie algebras
for the alloy $\frak g$ will be denoted by $\bold L(\frak g)$ and called
{\it the universal envelopping Lie algebra\/} of alloy $\frak g$.

\proclaim{Theorem 2} The representations of alloy $\frak g$ define
the representations of the uni\-ver\-sal envelopping Lie algebra $\bold L(\frak
g)$ and vice versa.
\endproclaim

\remark{Remark 10} The universal envelopping Lie algebras of alloys are
infinite-dimen\-sional as a rule.
\endremark

Note that the infinite-dimensional Lie algebras $\bold L(\frak g)$ are
essentially wild, however, the theory of their {\it finite-dimensional\/}
representations may be a bit timer (as a rule time infinite-dimensional
Lie algebras have only trivial in some sense finite-dimensional
representations). This statement will be illustrated by examples below.

Let $\frak g$ be an alloy and $\Cal L_0(\frak g)$ be the category of all
finite-dimensional envelopping Lie algebras for $\frak g$ (the elements of
$\Cal L_0(\frak g)$ are Lie algebras $\tilde\frak g$ with the fixed
imbedding $\iota:\frak g\mapsto\tilde\frak g$). Let us denote by
$\FR(\frak h)$ the category of finite-dimensional rep\-re\-sen\-tations of
a Lie algebra $\frak h$. If $\tilde\frak g_1,\tilde\frak g_2\in\Cal
L_0(\frak g)$ and $\pi_i\in\FR(\tilde\frak g_i)$ ($i=1,2$) then put
$\pi_1\underset\df\to\sim\pi_2$ iff $\left.\pi_1\right|_{\frak g}=
\left.\pi_2\right|_{\frak g}$.

\proclaim{Proposition 4} $\FR(\bold L(\frak g))$ is a tensor category,
which may be represented as
$$\FR(\bold L(\frak g))=\bigoplus_{\tilde\frak g\in\Cal L_0(\frak g)}
\FR(\tilde\frak g)/\sim.$$
\endproclaim

The proposition 4 means that the representations of different Lie algebras
$\tilde\frak g$ (supplied by with the fixed imbeddings of $\frak g$ into
them) allow mutual tensor pro\-ducts.

\subhead 1.5. Examples
\endsubhead
Let us consider a four-dimensional alloy $\asltwo$ generated by the elements
$e_0$, $f_{\pm}$ and $e_1$ with brackets $[e_0,f_{\pm}]=\pm f_{\pm}$,
$[f_+,f_-]=2e_0+e_1$. This is an alloy related to a cross-projective
representation of a pair $(\sltwo,\Bbb C)$ in general position (up to
an automorphism of $\sltwo$).

\proclaim{Theorem 3} Any irreducible finite-dimensional $e_0$-diagonal
representation of the alloy $\asltwo$ has the form
$$\aligned
f_+\mapsto\left(\matrix 0&A_1&\boldkey 0\\&\ddots&A_N\\\boldkey 0&&0
\endmatrix\right),&\quad
f_-\mapsto\left(\matrix 0&&\boldkey 0\\B_1&\ddots&\\\boldkey 0&B_N&0
\endmatrix\right),\quad\\
e_0\mapsto\left(\matrix C_0&&\boldkey 0\\&\ddots&\\\boldkey 0&&C_N
\endmatrix\right),&\quad
e_1\mapsto[f_+,f_-]-2e_0,
\endaligned$$
where $C_i=(\gamma-i)E$ ($E$ is the $n_i\times n_i$-dimensional unit matrix),
$A_i$, $B_i$ are $n_{i-1}\times n_i$- and $n_i\times n_{i-1}$-dimensional
matrices, respectively. The condition of irreducibility puts some relations
on the numbers $n_i$ (the dimensions of blocks), namely
$$n_0=n_N=1,\quad n_i\le n_{i-1}+n_{i+1}$$
as well as additional relations on matrices $A_i$ and $B_i$, namely that
the pair of $n_i\times n_i$ matrices $A_{i+1}B_{i+1}$ and $B_iA_i$
algebraically generate the whole matrix algebra $\Mat_{n_i}(\Bbb C)$ for
all $i$.
\endproclaim

\remark{Remark 11 (An exercise)} It is a nice exercise to classify the
irreducible rep\-re\-sen\-tations of $\asltwo$ of small dimensions and
to calculate the decompositions of their tensor products. My own
calculations gave me a lot of pleasure.
\endremark

\remark{Remark 12} It is important that the finite-dimensional irreducible
representations of $\asltwo$ have the continuous moduli. I suspect that
such situation is typical for representations of alloys and differs it from
the theory of finite-dimensional representations of finite-dimensional Lie
algebras.
\endremark

\Refs
\roster
\item"[1]" Juriev D., Topics in hidden symmetries. V: funct-an/9611003.
\item"[2]" Juriev D.V., Approximate representations and Virasoro algebra:
math.RT/9805001.
\item"[3]" Juriev D.V., $q_R$--conformal symmetries in two-dimensional nonlocal
quantum field theory, categorical representation theory and Virasoro algebra:
q-alg/9712009.
\item"[4]" Kolmogorov A.N., On a representation of continuous functions
of several variables as superpositions of continuous functions of the less
number of variables. Doklady AN SSSR. 108(2) (1956) 179-182.
\item"[5]" Arnold V.I., On a representation of continuous functions of three
variables as su\-per\-po\-sitions of continuous functions of two variables.
Metem.Sb. 48(1) (1959) 3-74, [E] 56(3) (1962) 382.
\endroster
\endRefs
\enddocument